\documentclass{siamltex}

\usepackage{amsmath, latexsym, amsfonts, amssymb, amscd,bbm}
\usepackage[english]{babel}
\usepackage[dvips]{graphicx}

\def \R{\mathbb R}

\def \P{\mathbb P}

\newcommand{\x}{\mathbf{x}}

\usepackage{color}

\newcommand{\e}{{\mathrm e}}

\newcommand{\norm}[1]{\left \| #1 \right \|}

\newcommand{\cN}{{\ensuremath{\mathcal N}} }

\newcommand{\M}{\mathfrak{M}}
\renewcommand{\H}{\mathfrak{B}}

\newcommand{\G}{{\mathcal G}}
\newcommand{\n}{{\mathbf n}}
\renewcommand{\x}{{\mathbf x}}
\newcommand{\X}{{\mathbf X}}

\newcommand{\james}[1]{{\color{black}#1}}

\title{Phase Reduction of Stochastic Biochemical Oscillators\thanks{PCB and JNM were supported
by the National Science Foundation (DMS-1613048).}}
\author{Paul C. Bressloff\thanks{Department of Mathematics, University of Utah, Salt Lake City, UT 84112 USA  ({\tt bressloff@math.utah.edu}).} \and James N. MacLaurin\thanks{Department of Mathematical Sciences, New Jersey Institute of Technology, Newark, NJ USA  ({\tt james.n.maclaurin@njit.edu}).}}

\begin{document}

\maketitle

\begin{abstract}

A common method for analyzing the effects of molecular noise in chemical reaction networks is to approximate the
 underlying chemical master equation by a Fokker-Planck equation, and to study the statistics of the associated
  chemical Langevin equation. This so-called diffusion approximation involves performing a perturbation expansion
   with respect to a small dimensionless parameter $\epsilon =\Omega^{-1}$, where $\Omega$ characterizes the
    system size. For example, $\Omega$ could be the mean number of proteins produced by a gene regulatory 
    network. In the 
  deterministic limit $\Omega\rightarrow \infty$, the chemical reaction network evolves according to a system of
   ordinary differential equations based on classical mass action kinetics.
   In this paper we develop a phase reduction method for chemical reaction networks that 
support a stable limit cycle in the deterministic limit. We present a variational principle for the phase reduction, yielding
 an exact analytic 
expression for the resulting phase dynamics. We demonstrate that this decomposition is accurate over 
timescales that are exponential in the system size $\Omega$. This contrasts with the phase equation obtained under 
the diffusion approximation, which is only accurate up to times $O(\Omega)$. In particular, we show that for a constant 
$C$,  
the \james{probability that the system leaves an $O(\zeta)$ neighborhood of the limit cycle before time $T$ scales
as $T\exp\big(-C\Omega b\zeta^2 \big)$, where $b$ is the rate of attraction to the limit cycle}. We illustrate our analysis using 
the example of a chemical Brusselator.

\end{abstract}

\renewcommand{\thefootnote}{\fnsymbol{footnote}}

\section{Introduction}

Genetically identical cells exposed to the same environmental conditions can show
significant variation in molecular content and marked differences in phenotypic characteristics.
This intrinsic variability is linked to to the fact that many cellular events at the genetic level involve small numbers of molecules (low copy numbers), and has led to a large number of studies focused on investigating the origins and consequences
of intrinsic or molecular noise in gene expression, see the reviews \cite{Collins05,Paulsson05,Raj08,Tsimiring14,Bressloff17}. 
In the presence of intrinsic noise, standard deterministic models of chemical reaction networks (CRNs), which are based on mass action kinetics, have to be replaced by discrete stochastic models. The stochastic dynamics of the CRN can then be formulated in terms of a chemical master equation, which determines the evolution of the probability distribution $P(\n,t)$ for $\n=(n_1,\ldots, n_K)$, where $n_i$ is the number of each molecular species labeled $i=1,\ldots,K$. In general, it is not possible to obtain exact solutions of the master equation, even in the case of a stationary solution. (Note, however, that recent progress has been made by generalizing the theory of deterministic chemical reaction networks to stochastic models, using notions of weak reversibility and deficiency zero \cite{Anderson10,Anderson15}.) Therefore, one often resorts to some form of approximation scheme. The most common approach is to carry out a Taylor expansion of the master equation with respect to a dimensionless variable $\Omega$ that characterizes the size of the system \cite{vanKampen92,Elf03}. This so-called diffusion approximation yields a Fokker-Planck (FP) equation, whose solution determines the probability density $p(\x,t)$ for the concentrations $\x=\n/\Omega$. 

A general issue is to what extent the system-size expansion of a chemical master equation accounts for the statistics of the underlying CRN. In cases where the deterministic mass action kinetics have a unique stable fixed point, the corresponding stationary distribution of the master equation can be well approximated by Gaussian-like fluctuations around the fixed point. One can thus obtain good estimates for the mean and variance of molecular concentrations using the diffusion approximation \cite{Swain02,Elf03,Paulsson05}. However, such an approximation can break down when there are multiple fixed points (metastability). Although one can still account for the  effects of fluctuations well within the basin of attraction of each metastable fixed point, there is now a small probability that there is a noise--induced transition to the basin of attraction of another fixed point. Since the probability of such a transition is usually of order $\e^{-\tau \Omega}$  with $\tau=O(1)$, except close to the boundary of the basin of attraction, such a contribution cannot be analyzed accurately using standard Fokker--Planck methods. That is, for large $\Omega$ the transitions between different metastable states typically involve rare transitions (large fluctuations) that lie in the tails of the associated probability distribution, where the Gaussian approximation no longer necessarily holds. These exponentially small transitions play a crucial role in allowing the system to approach the unique stationary state (if it exists) in the asymptotic limit $t\rightarrow \infty$. In other words, for bistable (or multistable) systems, the limits $t\rightarrow \infty$ and $N\rightarrow \infty$ do not necessarily commute \cite{Hanggi84,Vellela10}. In order to ensure accurate estimates of transition rates, one has to use alternative approximation schemes based on large deviation theory and perturbation methods, which are applied directly to the underlying master equation or CK equation \cite{Bressloff17}.
We also note that several recent papers bound the probability of large deviations over fixed timescales \cite{Kumar16,Agazzi17}.

Similar issues hold for other types of attractor such as stable limit cycles. In the absence of noise, one can represent the dynamics of a limit cycle in terms of motion around a closed curve in state space. Each point on the limit cycle can be assigned a phase $\theta\in [0,2\pi]$ such that the dynamics is given by
$\dot{\theta}=\omega_0$,
where $\omega_0=2\pi/\Delta_0 $ and $\Delta_0$ is the period of the oscillation. When noise is included, one tends to observe an irregular trajectory around the limit cycle. If the trajectory wanders too far from the limit cycle then it may escape completely, and the notion of phase breaks down. The mean time to escape a neighborhood of the limit cycle is a measure of robustness to transverse fluctuations. However, even if the trajectory remains close to the limit cycle, the resulting stochastic phase of the oscillator undergoes diffusion around the limit cycle, so all phase information is eventually lost. The effective diffusion coefficient of this process is one measure of robustness to phase decoherence \cite{Gonze02}. (For a general discussion of robustness in biological processes, see \james{\cite{Kitano04,Maini12}.}) A second effect of noise is that it can extend the parameter regime over which oscillations can occur. That is, even though the deterministic system converges to a fixed point, the system exhibits oscillatory behavior when noise is included. This can be established by looking at the power spectra of the concentrations \cite{Boland08,Boland09}. As far as we are aware, almost all previous studies of stochastic CRN oscillators are based on the diffusion approximation of chemical master equations or modifications of such approximations, which generate more accurate estimations of the statistics that remain uniformly accurate for long times \cite{Minas17}. One exception is a hybrid computational study, which combines simulations of the discrete Markov process via Gillespie and solutions of the approximate SDE \cite{Suvak12}. In addition, a recent paper \cite{keskin18} looks at noise in the vertebrate segmentation clock, employing a computational model to study biochemical oscillations in coupled cells.

In this paper, we focus on the robustness of limit cycle oscillators to transverse fluctuations, which we analyze by working directly with the underlying discrete Markov process rather than resorting to a diffusion approximation. The limitation of the latter within this context is that it is only accurate up to times $O(\Omega)$, whereas we are interested in obtaining exponential bounds on the typical time that the system spends in a neighborhood of the limit cycle.
Indeed, we prove that, with very high probability, this time scales as $\exp\big( C\Omega b \big)$, where $C$ is a constant, and $b$ is the rate of attraction towards the limit cycle. The fact that the timescale of the bound is so large demonstrates the utility of using an exact representation of the phase.  We begin in \S 2 by defining a CRN and reviewing the method of isochrons for analyzing deterministic biochemical oscillators. We then use the method of isochrons to derive a stochastic phase equation under the diffusion approximation. In \S 3 we show how to derive a corresponding phase equation for the exact discrete Markov process using Kurtz's random time change representation \cite{Kurtz71,Kurtz76,Kurtz80,Anderson11}. This analysis is placed on a more rigorous footing in \S 4 by extending a variational principle that we recently developed for stochastic differential equations (SDEs) \cite{Bressloff18a}. The latter generates a unique amplitude-phase decomposition of a stochastic oscillator by minimizing the distance of the stochastic trajectory from the limit cycle with respect to an appropriately weighted norm. Applying the variational approach to the stochastic CRN model, we show that (i) sufficiently close to the limit cycle, the phase corresponds to the isochronal phase, and (ii) one can obtain rigorous exponential bounds on the probability of escape. (The latter is carried out in \S 5). We illustrate our theory by numerically simulating the stochastic Brusselator. Finally, in \S 6, we summarize the paper and draw some conclusions.

\section{Biochemical oscillators and the diffusion approximation}

\subsection{Chemical reaction network} Consider a well-mixed compartment containing a set of chemical species $X_j$, $j=1,\ldots,K$. Let $n_j$ be the number of molecules of $X_j$ and set ${\bf n}=(n_1,\ldots,n_K)$. A typical single-step chemical reaction takes the form
\begin{equation*}
s_1X_1+s_2X_2+\ldots \rightarrow r_1 X_1+r_2X_2+\ldots,\end{equation*}
where $s_j,r_j$ are known as {\em stoichiometric coefficients}. When one such reaction occurs the state ${\bf n}$ is changed according to
\[n_i\rightarrow n_i+r_i-s_i.\]
Introducing the vector ${\bf S}$ with $S_i=r_i-s_i$, we have ${\bf n}\rightarrow {\bf n}+{\bf S}$. The reverse reaction
\[\sum_jr_jX_j\rightarrow \sum_js_fX_j,\]
would then have ${\bf n}\rightarrow {\bf n}-{\bf S}$. More complicated multi-step reactions can always be decomposed into these fundamental single-step reactions with appropriate stoichiometric coefficients. 

In the case of a large number of molecules, one can describe the dynamics of a single step chemical reaction in terms of a kinetic or rate equation involving the concentrations $x_j=n_j/\Omega$ -- the {\em law of mass action}. Here $\Omega$ is a dimensionless quantity representing the system size, which is usually taken to be the mean number of molecules or some volume scale factor. In a well-mixed container there is a spatially uniform distribution of each type of molecule, and the probability of a collision depends on the probability that each of the reactants is in the same local region of space. Ignoring any statistical correlations, the latter is given by the product of the individual concentrations. It then follows that the kinetic equations take the form
\begin{equation}
\label{ostoch}
\frac{dx_i}{dt}=\kappa (r_i-s_i)\prod_{j=1}^Kx_j^{s_j}\equiv S_i \lambda(\x),
\end{equation}
where $\x=(x_1,\ldots,x_K)$ and $\kappa$ is a rate constant. Now suppose that there are $a=1,\ldots, M$ separate single-step reactions. Then
\begin{equation}
\label{rateg}
\frac{dx_i}{dt}:=F_i(\x)=\sum_{a=1}^MS_{ia}\lambda_a(\x),\quad i=1,\dots,K
\end{equation}
where $a$ labels a single-step reaction and ${\bf S}$ is the so-called $K\times M$ stoichiometric matrix for $K$ molecular species and $R$ reactions. Thus $S_{ia}$ specifies the change in the number of molecules of species $i$ in a given reaction $a$. The functions $\lambda_a$ are known as transition intensities or {\em propensities}. 

Mathematically speaking, the kinetic equations hold in the thermodynamic limit $\Omega \rightarrow \infty$. For finite $\Omega$, it is necessary to take into account fluctuations in the number $N_i(t)$ of each chemical species. Let $P(\n,t)=\P[N_1(t)=n_1,\ldots,N_K(t)=n_K]$ with $\n=(n_1,\ldots,n_K)$. Given the kinetic equations (\ref{0chemg}), the probability distribution $P(\n,t)$ evolves according to the chemical master equation takes the form
\begin{equation}
\label{masterg}
\frac{dP(\n,t)}{dt}=\Omega \sum_{a=1}^R\left (\prod_{i=1}^K{\mathbb E}^{-S_{ia}}-1\right )\lambda_a(\n/\Omega)P(\n,t),
\end{equation}
Here ${\mathbb E}^{-S_{ia}}$ is a step or ladder operator such that for any function $g(\n)$,
\begin{equation}
{\mathbb E}^{-S_{ia}}g(n_1,\ldots,n_i,\ldots,n_K)=g(n_1,\ldots, n_i-S_{ia},\ldots,n_K).
\end{equation}
One point to note is that when the number of molecules is sufficiently small, the characteristic form of a propensity function $\lambda_a(\x)$ in equation (\ref{masterg}) has to be modified:
\[ \left (\frac{n_j}{\Omega}\right )^{s_j }\rightarrow \frac{1}{\Omega^{s_j}}\frac{n_j!}{(n_j-s_j)!}.\]
In general, it is not possible to obtain exact solutions of the master equation (\ref{masterg}) even in the case of a stationary solution. (Note, however, that recent progress has been made by generalizing the theory of deterministic chemical reaction networks to stochastic models, using notions of weak reversibility and deficiency zero \cite{Anderson10}.) Therefore, one often resorts to some form of approximation scheme.

\subsection{Biochemical oscillators and isochrons} The next step is to assume that the deterministic dynamical system (\ref{rateg}) supports a stable limit cycle. That is, there exists a stable periodic solution ${\bf x}=\Phi(t)$ with $\Phi(t)=\Phi(t+\Delta_0)$, where
$\omega_0=2\pi/\Delta_0$ is the natural frequency of the oscillator. The dynamics on the limit cycle can be described by a uniformly rotating phase such that
\begin{equation}
\frac{d\theta}{dt}=\omega_0,
\end{equation}
and $x={\Phi}(\theta(t))$ with ${\Phi}$ a $2\pi$-periodic function and $\theta(t)=\theta_0+\omega_0 t$. The dynamical equations for $\Phi$ can thus be written as
\begin{equation}
\label{fi}
\omega_0 \frac{d\Phi}{d\theta} = {F}(\Phi(\theta)).\end{equation}
Differentiating both sides with respect to $\theta$ gives
\begin{equation}
\frac{d}{d\theta} \left (\frac{d\Phi}{d\theta}\right )=\omega_0^{-1}{J}( \theta)\cdot \frac{d\Phi}{d\theta},
\label{nadj}
\end{equation}
where
${J}$ is the $2\pi$-periodic Jacobian matrix 
\begin{equation}
{J}_{jk}(\theta)\equiv \left . \frac{\partial{F}_j}{\partial x_k}\right |_{x=\Phi(\theta)}.
\end{equation}
Since the dynamical system is autonomous, it follows that it is neutrally stable with respect to perturbations along the limit cycle.

Our main concern in this paper is to to characterize how intrinsic noise associated with fluctuations in the number of each chemical species for finite $\Omega$ affects the phase of the oscillator. This will require understanding how the deterministic system responds to small perturbations. Therefore, consider the perturbed deterministic dynamical system
\begin{equation}
\label{rateg2}
\frac{dx_i}{dt}:=F_i(\x)+\epsilon G_i(\x),\quad i=1,\dots,K,
\end{equation}
where $0<\epsilon \ll 1$ and $G_i(\x)$ is bounded. A classical method for derived the corresponding phase equation is based on the construction of isochrons \cite{Winfree80,Kuramoto84,Glass88,Holmes04,Ashwin16}.
Suppose that the unperturbed system (\ref{rateg2}) is observed stroboscopically at time intervals of length $\Delta_0$. This leads to a
Poincare mapping
\begin{equation*}
\x(t)\rightarrow \x(t+\Delta_0)\equiv {\mathcal P}(\x(t)),
\end{equation*}
which has all points on the limit cycle as fixed points. The iscochron through a point $\x^*$ on the limit cycle is the $(K-1)$-dimensional hypersurface ${\mathcal I}$, consisting of
all points in the vicinity of $\x^*$ that are attracted to $\x^*$ under the action of ${\mathcal P}$. A unique isochron can be drawn through each point on the limit cycle (at least locally) so the isochrons can be parameterized by the phase,
${\mathcal I}={\mathcal I}(\theta)$. Finally, the definition of phase is extended by taking all points $\x\in {\mathcal I}(\theta)$ to have the same phase,
$\Theta(\x)=\theta$, which then rotates at the natural frequency $\omega_0$ (in the
unperturbed case). Hence, for an unperturbed oscillator in the vicinity of the limit cycle we have
\begin{eqnarray}
\omega_0 = \frac{d\Theta}{dt}=\sum_{j=1}^K\frac{\partial \Theta}{\partial x_j}\frac{dx_j}{dt} =\sum_{j=1}^K\frac{\partial \Theta}{\partial x_j}F_j(\x) .\nonumber
\end{eqnarray}

Extending the definition of isochrons to the perturbed system (\ref{rateg2}), we have to leading order
\begin{equation}
\frac{d\Theta}{dt}=\sum_{j=1}^K\frac{\partial \Theta}{\partial x_j}(F_j(\x)+\sqrt{\epsilon} G_j({\x},t))=\omega_0+\sqrt{\epsilon}\sum_{j=1}^K\ \frac{\partial \Theta}{\partial x_j}G_j(\x,t). \nonumber
\end{equation}
As a further leading order approximation, deviations of $\x$ from the limit cycle are ignored on the right-hand side. Hence, setting $\x(t)=\Phi(\theta(t))$ with $\Phi$ the $2\pi$-periodic solution on the limit cycle,
\begin{equation}
\frac{d\Theta}{dt}=\omega_0+\sqrt{\epsilon}\sum_{j=1}^K \left .\frac{\partial \Theta}{\partial x_k}\right |_{\x=\Phi}G_k(\Phi,t) .\nonumber
\end{equation}
Finally, since points on the limit cycle are in 1:1 correspondence with the phase $\theta$, one can set $\Theta(\Phi(\theta))=\theta$ to obtain the closed
phase equation
\begin{equation}
\frac{d\theta}{dt}=\omega_0+\sqrt{\epsilon} \sum_{j=1}^K R_j(\theta)G_j(\Phi(\theta),t)
\label{cbphase}
\end{equation}
where
\begin{equation}
R_k(\theta)= \left . \frac{\partial \Theta}{\partial x_k}\right |_{\x=\Phi(\theta)}
\label{Q2}
\end{equation}
is a $2\pi$-periodic function of $\theta$ known as the $k$th component of the {\em phase response curve} (PRC).
One way to calculate the PRC $R(\theta)$ is to note that it is a $2\pi$-periodic solution of the linear equation \cite{Erm96,Erm10,nakao2016}
\begin{equation}
\label{adj}
\omega_0\frac{dR(\theta)}{d\theta}=-J(\theta)^{\top}\cdot  R(\theta),
\end{equation}
with the normalization condition
\begin{equation}
R(\theta)\cdot \frac{d\Phi(\theta)}{d\theta}=1.
\end{equation}
Here $J(\theta)^{\top}$ is the transpose of the Jacobian matrix $J(\theta)$, i.e.
\begin{equation}
\label{Jac}
J_{jk}(\theta)\equiv \left . \frac{\partial F_j}{\partial x_k}\right |_{\x=\Phi(\theta)}.
\end{equation}

 \subsection{Diffusion approximation and phase reduction}
 
\begin{figure}[t!]
\centering
\includegraphics[width=8cm]{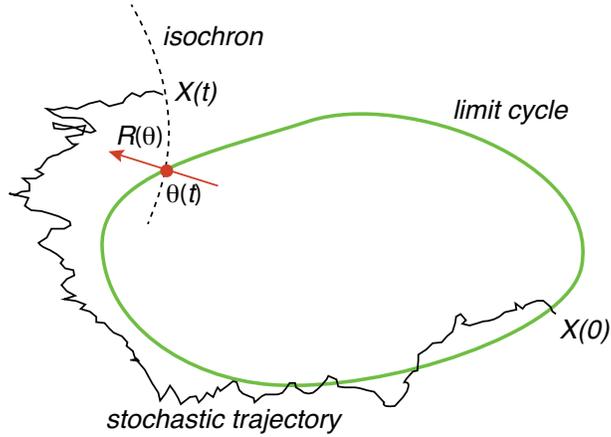}
\caption{Stochastic trajectory of chemical reaction network for finite $\Omega$, assuming that the corresponding mass action kinetics supports a stable limit cycle. The method of isochrons determines the phase $\theta(t)$ by tracing where the isochron through $X(t)$ intersects the limit cycle.The response to perturbations depends on the phase response curve ${R}(\theta)$, which is normal to the isochron at the point of intersection with the limit cycle.}
\label{scycle}
\end{figure}

 Now suppose that we include the effects of intrinsic noise by taking $\Omega$ to be large but finite, with the corresponding deterministic system operating in an oscillatory regime. Intuitively, one might expect that simulating the chemical reaction network using Gillespie's stochastic simulation algorithm \cite{Gillespie77,Gillespie01,Gillespie13}, for example, will generate a stochastic trajectory that remains in a neighborhood of the limit cycle (up to some stopping time), as illustrated schematically in Fig. \ref{scycle}. Given the marginal stability of the phase in the deterministic system, it is natural to decompose the effects of noise into longitudinal (phase) and transverse (amplitude) fluctuations of the limit cycle \cite{Gonze02,Boland09,Koeppl11,Bonnin17,Bressloff18a}. The basic intuition is that Gaussian-like transverse fluctuations are distributed in a tube of radius $1/\sqrt{\Omega}$ whereas the phase around the limit cycle undergoes undamped Brownian motion. Thus, the stochastic vector $\X(t)={\bf N}(t)/\Omega$ is decomposed in the form\footnote{Note that $N_j(t)$, $j=1,\ldots,K$ are discrete random variables so that $X_j(t)$ is also discrete. However, if $\Omega$ is sufficiently large, then the trajectory consists of small jumps so that it can be approximated by a continuous trajectory, assuming that the number of jumps in a small time interval is bounded appropriately.}
\begin{equation}
\label{defcon}
\X(t)=\Phi(\theta(t))+\frac{1}{\sqrt{\Omega}}{\bf v}(t),
\end{equation}
where the scalar random variable $\theta(t)-\omega_0 t$ represents the undamped random phase shift along the limit cycle, and ${\bf v}(t)$ is a transversal perturbation, see Fig. \ref{scycle}. However, it is important to note that the decomposition (\ref{defcon}) is not unique, so that the precise definition of the phase depends on the particular method of analysis. One possibility is to identify $\theta(t)$ with the isochronal phase, as shown in Fig. \ref{scycle}. An alternative approach is to introduce a weighted inner product on $\R^K$ and to specify ${\bf v}(t)$ by projecting the full solution on to the limit cycle using Floquet vectors. In a sufficiently small neighborhood of the limit cycle, the resulting phase variable coincides with the isochronal phase \cite{Bonnin17,Bressloff18a}. This has the advantage that the amplitude and phase decouple to linear order in $1/\Omega$.

As we noted in the introduction, most previous studies of stochastic CRN oscillators are based on the diffusion approximation of chemical master equations \cite{Gonze02,Boland08,Boland09,Minas17}. The concentration vector $\X(t)$ then evolves according to a stochastic differential equation (SDE). It is possible to derive the SDE rigorously using Kurtz's random time change representation \cite{Kurtz71,Kurtz76,Kurtz80,Anderson11}. A more heuristic approach is to perform a system-size expansion of the master equation \cite{Gardiner09,vanKampen92,Elf03}.
That is, set $\lambda_a(\n/\Omega)P(\n,t)\rightarrow \lambda_a(\x)p(\x,t)$ with $\x=\n/\Omega$, and treat $\x$ as a continuous vector. Then
\begin{eqnarray*}
 \prod_{i=1}^K{\mathbb E}^{-S_{ia}}h(\x)&=&h(\x-{\bf S}_a/\Omega)\\
 &=&h(\x)-\Omega^{-1}\sum_{i=1}^KS_{ia}\frac{\partial h}{\partial x_i}+\frac{1}{2\Omega^2}\sum_{i,j=1}^KS_{ia}S_{ja}\frac{\partial^2h(\x)}{\partial x_i\partial x_j}+O(\Omega^{-3})
\end{eqnarray*}
Carrying out a Taylor expansion of the master equation to second order thus yields the multivariate Fokker-Planck (FP) equation
\begin{equation}
\label{FPg}
\frac{\partial p}{\partial t}=-\sum_{i=1}^K\frac{\partial F_i(\x)p(\x,t)}{\partial x_i}+\frac{1}{2\Omega}\sum_{i,j=1}^K\frac{\partial^2 D_{ij}(\x)p(\x,t)}{\partial x_i\partial x_j},
\end{equation}
where
\begin{equation}
\label{VBg}
F_i(\x)=\sum_{a=1}^R S_{ia}\lambda_a(\x),\quad D_{ij}(\x)=\sum_{a=1}^R S_{ia}S_{ja}\lambda_a(\x).
\end{equation}
The FP equation (\ref{FPg}) corresponds to the multivariate Ito stochastic differential equation (SDE)
\begin{equation}
\label{Langg}
dX_i=F_i(\X)dt+\frac{1}{\sqrt{\Omega}}\sum_{a=1}^R B_{ia}(\X)dW_a(t),
\end{equation}
where $W_a(t)$ are independent Wiener processes \cite{Gardiner09},
\begin{equation}
\langle dW_a(t)\rangle = 0,\quad \langle dW_a(t)dW_b(t')\rangle =\delta_{a,b}\delta(t-t')dt \, dt',\end{equation}
 and ${\bf D}={\bf B}{\bf B}^T$, that is,
\begin{equation}
B_{ia}=S_{ia}\sqrt{\lambda_a(\x)}.
\end{equation}
One can now carry out an amplitude-phase decomposition of the SDE (\ref{Langg}) to estimate the size of Gaussian fluctuations in the longitudinal and transverse directions \cite{Gonze02,Boland08,Boland09}. It is also possible to modify the diffusion-based approximation, in order to obtain more accurate estimations of the statistics, which remain uniformly accurate for long times \cite{Minas17}.

In terms of the isochronal phase reduction method highlighted above, one could simply treat the sum over Wiener processes in equation (\ref{Langg}) as a perturbation of the deterministic equation, by setting $\epsilon = 1/\Omega$ and
\[G_j(\x)dt=\sum_{a=1}^R B_{ia}(\X)dW_a(t).\]
This would then lead to a stochastic version of the phase equation (\ref{cbphase}). However, one needs to be careful in applying the phase reduction procedure to the Ito SDE. First, in order to apply the usual rules of calculus it is necessary to convert the Ito SDE to Stratonovich form using Ito's lemma. This would then yield an $O(1/\Omega)$ correction to the drift term. Since we are interested in the leading order contributions to the phase equation, we can ignore such a term, and simply take the phase to evolve according to the Ito SDE
\begin{equation}
 d\theta =\omega_0dt+\sqrt{\frac{1}{\Omega}} \sum_{j=1}^K R_j(\theta)\sum_{a=1}^R B_{ja}(\Phi(\theta(t)))dW_a(t).
\label{phaseFP}
\end{equation}
The second, and more significant issue, is that even for large $\Omega$ there is a small but non-zero probability that a stochastic trajectory generated by the SDE (\ref{Langg}) leaves any neighborhood of the limit cycle, in which case the notion of phase breaks down. In previous work, we have developed a variational principle for SDEs that provides strong exponential bounds on the escape probability \cite{Bressloff18a}. (We have also extended the variational formulation to the case of piecewise-deterministic Markov processes \cite{Bressloff18b}.) One could directly apply our analysis to the SDE (\ref{Langg}). However, the latter is already an approximation of the exact discrete Markov process. This motivates the analysis of this paper, namely, to derive an explicit phase equation to $O(1/\Omega)$, and to obtain bounds on the probability of escape by working directly with the discrete Markov process, rather than its diffusion approximation. The mathematical framework underlying our approach is the random time change representation of chemical reaction networks, which has been developed by Kurtz \cite{Kurtz71,Kurtz76,Kurtz80,Anderson11}.
As we explained in the introduction, a fundamental aim in systems biology is to understand the robustness of the oscillation: that is, for a given system size, how long should one typically expect the oscillation to persist? The diffusion approximation is limited because it only aims to reproduce the statistics of the oscillation over timescales of $O(\Omega)$. However, as we show in the following sections, one should expect the oscillation to persist for exponentially long periods of time; a metastable timescale over which the diffusion approximation breaks down. Moreover, our variational definition of the phase will be accurate for the entire time that the system is in a neighborhood of the limit cycle.

\section{Random time change representation and phase reduction}
Let the random variable $N_i(t)$ denote the number of molecules in species $i$. An explicit representation of $N_i(t)$ is given by
\begin{equation}
\label{rtc0}
N_i(t) = N_i(0) + \sum_{a=1}^MS_{ia}\cN_a(t),
\end{equation}
where $\cN_a(t)$ denotes the number of times reaction $a$ has occurred by time $t$, and is a jump process for which the jump rate is locally given by $\Omega f_a({\bf N}(t)/\Omega)$. {It can be proved that if $Y(t)$ is a unit rate Poisson process, that is,
\[\P[Y(t)=m|Y(0)=0] = \frac{t^n\e^{- t}}{n!},\]
then a representation of $\cN_a(t)$ is given by the \emph{time change representation} \cite{Kurtz71,Anderson11},
\begin{equation*}
\cN_a(t) = Y_a\left(\Omega\int_0^t\lambda_a({\bf N}(s)/\Omega)ds \right). 
\end{equation*}
It follows that the expected number of reactions will grow linearly with $\Omega$. Moreover, using the law of large numbers for Poisson processes, 
\begin{equation*}
   \lim_{\Omega \to\infty}\frac{Y_a(\Omega T) }{\Omega}= T,
\end{equation*}
one has
\begin{equation*}
\lim_{\Omega \to \infty}\frac{ \cN_a(t)}{\Omega}  =  \int_0^t\lambda_a(\x(s))ds
\end{equation*}
with $\x(s)=\lim_{\Omega \rightarrow \infty}{\bf N}(s)/\Omega$. This fact can be used to derive the deterministic equations (\ref{0chemg}) in the thermodynamic limit. 

As before, we define $X_i(t) =N_i(t)/\Omega$ to be the concentration of species $i$ at time $t$, which has the deterministic dynamics (\ref{rateg}) in the large $\Omega$ limit. Let ${\bf S}_a = \lbrace S_{ia} \rbrace_{i=1,\ldots, S}$ be the vector representing the change in species number by reaction $a$, and $\grave{Y}_a(z) = Y_a(z) - z$ be the $a^{th}$ Poisson process minus its mean. We then find that
\begin{align}\label{rtc}
\X(t) = \X(0) + \int_0^t {\bf F}\big(\X(s) \big) ds + \Omega^{-1}\sum_{a=1}^M \grave{Y}_a\bigg(\Omega \int_0^t \lambda_a\big( \X(s) \big)ds \bigg){\bf S}_a.
\end{align}
Hence, we have represented the evolution of the concentration as the sum of a deterministic dynamics, and a stochastic noise of zero mean. (Taking the stochastic term to have zero mean will be useful later, when we use the theory of martingales, see \S 5.) It has been proved that in the large $\Omega$ limit, the dynamics converges to the ODE \cite{Kurtz71}
\begin{equation}\label{0chemg}
\X(t) = \X(0) + \int_0^t {\bf F}\big(\X(s) \big) ds.
\end{equation}

We now sketch a heuristic method for deriving a phase equation from (\ref{rtc}). In \S 4 we will put our derivation on a more rigorous footing using a variational approach. One way to specify the Poisson process $\cN_a(t) $ is in terms of the event or jump times $t^a_k$, $k=1,2,\ldots$. That is,
\begin{equation}
\label{Ra}
\cN_a(t)=\sum_{k\geq 1}H(t-t^a_k),\end{equation}
where the inter event times $t^a_{k+1}-t^a_k$ are generated from an exponential distribution:
\begin{equation*}
\P(\tau \leq t^a_{k+1}-t^a_k \leq \tau + \Delta \tau|t^a_k)=\lambda_a(\X(\tau)) \Delta \tau \, \e^{-\Lambda(t^a_k,\tau)},
\end{equation*}
with
\[\Lambda(t,\tau)=\int_t^{t+\tau}\lambda_a(\X(s))ds.\]
Formally speaking, $dH(t-t^a_k)=\delta(t-t^a_k)dt$. Differentiating equation (\ref{rtc}) thus gives
\begin{equation}
\frac{d\cN_a}{dt}=\sum_{k\geq 1}\delta(t-t^a_k).
\end{equation}
\begin{equation}
\frac{dX_i}{dt}=F_i(\X)+\Omega^{-1}\sum_{a=1}^M {S}_{ia}\left (\sum_k\delta(t-t^a_k)-\Omega \lambda_a(\X(t))\right ).
\end{equation}

If we treat the sum on the right-hand side of the above equation as a perturbation of the deterministic equation, then we can formally use the isochronal phase reduction method to derive the following stochastic phase equation to $O(1/\Omega)$:
\begin{equation}
\frac{d\theta}{dt}=\omega_0+{\frac{1}{\Omega}} \sum_{j=1}^K R_j(\theta(t))\sum_{a=1}^M {S}_{ja}\left (\sum_k\delta(t-\overline{t}^a_k)-\Omega \lambda_a(\Phi(\theta(t)))\right ).
\label{phaseP}
\end{equation}
Here $\overline{t}_k^a$ are the sequence of jump times for the inhomogeneous Poisson process 
\begin{equation}
\label{Ya}
\overline{Y}_a(t)=Y_a\bigg(\Omega \int_0^t \lambda_a\big( \Phi(\theta(s)) \big)ds \bigg).\end{equation}
(In general, the times $\overline{t}_k^a$ will differ from the exact times ${t}_k^a$, since in the former case $X$ is evaluated on the limit cycle.) In terms of differentials, \james{(see \cite{Protter05} or $\S$A.3 of \cite{Anderson15} for a definition of the stochastic integral with respect to $\overline{Y}$)}
\begin{equation}
d\theta(t)=\omega_0 t+{\frac{1}{\Omega}} \sum_{j=1}^K R_j(\theta(t))\sum_{a=1}^M {S}_{ja}\left \lbrace d\overline{Y}_a(t)-\Omega \lambda_a(\Phi(\theta(t)))dt\right \rbrace.
\end{equation}
Finally, integrating equation (\ref{phaseP}) with respect to $t$, we obtain the stochastic phase equation
\begin{equation}
\label{theta0}
\theta(t)=\theta(0)+\omega_0t +{\frac{1}{\Omega}} \sum_{j=1}^K\sum_{a=1}^M {S}_{ja} \int_0^tR_j(\theta(s))\left (d\overline{Y}_a(s)-\Omega \lambda_a(\Phi(\theta(s)))ds\right ),
\end{equation}
or, equivalently
\begin{equation}
\label{theta}
\theta(t)=\theta(0)+\omega_0t +{\frac{1}{\Omega}} \sum_{j=1}^K\sum_{a=1}^M {S}_{ja} \left [\sum_{k\geq 1} R_j(\theta(\overline{t}_a^k))H(t-\overline{t}_a^k)-\Omega \int_0^t R_j(\theta(s))\lambda_a(\Phi(\theta(s)))ds\right ]. \end{equation}

\subsection{Bounding the number of reactions over a time interval}
We expect the above phase equation to be accurate up to times $O(\Omega)$ provided there is an exponentially small probability that the number of reactions is arbitrarily large. This in turn requires bounds on the reaction rates $\lambda_a$ in some neighborhood ${\mathcal U}$ of the limit cycle, which could be identified with the basin of attraction. (Once the system leaves the neighborhood of the limit cycle, one cannot sensibly define a phase dynamics in the first place.
The analysis in subsequent sections will demonstrate that the system stays in a neighborhood of the limit cycle for exponentially long times, with very high probability.) Therefore, 
we assume the following uniform bounds on the reaction rates $\lambda_a$,
\begin{align}
\sup_{1\leq a \leq M}\sup_{\theta \in [0,2\pi]}\sup_{z : \norm{z - \Phi(\theta)}_{\theta} \leq \eta}\big| \lambda_a\big(z \big) \big| \leq \overline{\lambda} \label{assumption lambda max}
\end{align}
\james{Here $\norm{\cdot}_{\theta}$ is a weighted norm that we define in the next section and $\eta> 0$ is a constant denoting the boundary of the attracting manifold of the limit cycle.}
Let $\mathcal{M}^a_{u,t}$ be the number of times that reaction $a$ occurs over a time interval $[u,t]$; an estimate of the error in equation (\ref{theta}) is $O({\mathcal M}_{0,t}/\Omega^2)$, where ${\mathcal M}_{0,t}=\sum_{a=1}^M{\mathcal M}^a_{0,t}$. In fact, for the error bounds in the coming sections, we need to bound the number of reactions that occur over arbitrary intervals of a certain length, i.e. we need a bound for arbitrary $u \in [0,t]$. To this end, for a constant integer $c>0$,
\begin{align*}
\mathbb{P}\bigg( \mathcal{M}_{u,t}^a \geq c \bigg) = \mathbb{P}\bigg( Y_a\big(\Omega\int_0^t \lambda_a(\X(s)))ds \big)-Y_a\big(\Omega\int_0^u \lambda_a(\X(s))ds \big) \geq c \bigg) \\
\leq \mathbb{P}\bigg( Y_a\big(\Omega\int_0^u \lambda_a(\X(s))ds+\Omega\overline{\lambda}(t-u) \big)-Y_a\big(\Omega\int_0^u \lambda_a(\X(s))ds \big) \geq c \bigg) \\
= \frac{1}{c!}\big\lbrace \Omega\overline{\lambda}(t-u) \big\rbrace^c \exp\bigg(-\Omega\overline{\lambda}(t-u) \bigg),
\end{align*}
using the definition of the Poisson distribution.  We find that for $c > \Omega\overline{\lambda}_a(t-u)$,
\begin{equation}\label{eq: bound number of reactions}
\mathbb{P}\big( \mathcal{M}_{u,t}^a \geq c \big) \leq \exp\big(-\widetilde{c}\Omega(t-u)^{-1} \big),
\end{equation}
for another constant $\widetilde{c} > 0$ which depends on $c$. Note that we have not explicitly written how $\widetilde{c}$ depends on $c$; it is sufficient for our purposes that the above relationship is exponential.

\section{Amplitude-phase decomposition}
In order that we may easily bound the probability of the stochastic CRN staying near the limit cycle oscillator, it is necessary to consider the amplitude-phase decomposition (\ref{defcon}) and define an appropriately weighted norm \james{$\|\X(t)-\Phi(\theta(t))\|_{\rho}$}, which is a measure of the distance of the stochastic trajectory from the limit cycle. However, as we have previously mentioned, the amplitude-phase decomposition is not unique, since it depends on the particular definition of the phase $\theta(t)$. Recently, we introduced a variational principle for SDEs that uniquely specifies the amplitude-phase decomposition by minimizing \james{$\|\X(t)-\Phi(\theta(t))\|_{\rho}$} for a specific choice of $\rho$, which has two significant features: (i) sufficiently close to the limit cycle, the phase corresponds to the isochronal phase, and (ii) one can obtain rigorous exponential bounds on the probability of escape \cite{Bressloff18a,Bressloff18b}. We will apply the variational approach to the stochastic CRN model, and show that both of these results carry over; (ii) will be established in \S 5. In order to specify the weight $\rho$, we first need to introduce a little Floquet theory.

\subsection{Floquet decomposition}
For any $0 \leq t$, define $\Pi(t) \in \R^{d\times d}$ to be the following Fundamental matrix for the ODE
\begin{equation}\label{eq: ODE Jz}
\frac{d z}{dt} = J(t)z
\end{equation}
for $J(t)=J (\omega_0 t)$. That is, $\Pi(t):=  \big( z_1(t) | z_2(t) | \ldots |z_d(t) \big)$, where $z_i(t)$ satisfies \eqref{eq: ODE Jz}, $z_1(0) = \Phi'(0)$, and $\lbrace z_i(0) \rbrace_{i=1}^d$ is an orthogonal basis for $\R^d$. Floquet Theory states that there exists a diagonal matrix $\mathcal{S}=\mbox{diag}(\nu_1,\ldots,\nu_d)$ whose diagonal entries are the Floquet characteristic exponents, such that
\begin{equation}
\label{pip}
\Pi\big( t\big) = P\big(\omega_0 t \big)\exp\big(t\mathcal{S} \big)P^{-1}(0),
\end{equation}
with $P(\theta)$ a $2\pi$-periodic matrix whose first column is $\Phi'(\omega_0t)$, and $\nu_1 = 0$. That is, $P(\theta)^{-1}\Phi'(\theta) ={\bf e}$ with ${\bf e}_{j}=\delta_{1,j}$. In order to simplify the following notation, we will assume throughout this paper that the Floquet multipliers are real and hence $P(\theta)$ is a real matrix. One could readily generalize these results to the case that $\mathcal{S}$ is complex. The limit cycle is taken to be stable, meaning that for a constant $b > 0$, for all $2\leq i \leq d$,
\begin{equation}\label{eq: nu bound}
\nu_i \leq - b.
\end{equation}
Since $F \in \mathcal{C}^2$, it must be that $P \in \mathcal{C}^2$. Furthermore $P^{-1}(\theta)$ exists for all $\theta$, since $\Pi^{-1}(t)$ exists for all $t$.

The above Floquet decomposition motivates the following weighted inner product: For any $\theta \in \R$, denoting the standard Euclidean dot product on $\mathbb{R}^d$ by $\langle \cdot , \cdot\rangle$,
\[
\langle u,v \rangle_\theta = \big\langle P^{-1}(\theta)u,P^{-1}(\theta)v \big\rangle,
\]
and $\norm{u}_\theta = \sqrt{\langle u,u\rangle_{\theta}}$. This weighting is useful for two reasons: it leads to a leading order separation of the phase from the amplitude, and it facilitates the strong bounds derived below, because the weighted amplitude always decays, no matter what the phase. 
The former is a consequence of the fact that the matrix $P^{-1}(\theta)$ generates a coordination transformation in which the phase in a neighborhood of the limit cycle coincides with the asymptotic phase defined using isochrons (see also \cite{Bonnin17}). This is reflected by the following relationship between the tangent vector to the limit cycle, $\Phi'(\theta)$, and
 the PRC $R(\theta)$ of equation (\ref{Q2}):
 \begin{equation}
\label{Rtan}
P^{\top}(\theta){R}(\theta)=P^{-1}(\theta)\Phi'(\theta),
\end{equation}
where
\begin{equation}
\H(\theta):= \norm{P^{-1}(\theta)\Phi'(\theta)}^2.
\end{equation}

We will proceed by defining $R(\theta)$ according to equation (\ref{Rtan}) and showing that it satisfies the adjoint equation (\ref{adj}). We will need the relation
\begin{equation}\label{eq: derivative P 0}
\omega_0 P'(\theta) = J(\theta)P(\theta)- P(\theta)\mathcal{S},
\end{equation}
which can be obtained by differentiating \eqref{pip}.
Differentiating both sides of equation (\ref{Rtan}) with respect to $\theta$, we have
\begin{align}
\label{Rtan2}
P^TR'+(P^{\top})'R=P^{-1}\Phi''+(P^{-1})'\Phi'.
\end{align}
Equation (\ref{eq: derivative P 0}) implies that
\begin{align*}
\omega_0 (P^{\top}(\theta))' = P^{\top}(\theta)J^{\top}(\theta)- \mathcal{S}P^{\top}(\theta)
\end{align*}
and
\[\omega_0(P^{-1}(\theta))'=-P^{-1}(\theta)J(\theta)+\mathcal{S}P^{-1}(\theta).\]
We have used the fact that ${\mathcal S}$ is a diagonal matrix and $P^{-1}P'+(P^{-1})'P=0$ for any square matrix. Substituting these identities in equation (\ref{Rtan2}) yields
\begin{align*}
& P^T(R'+\omega_0^{-1}J^{\top}R)-\omega_0^{-1} \mathcal{S}P^{\top}R =P^{-1}[\Phi''-\omega_0^{-1}J\Phi']+\omega_0^{-1}\mathcal{S}P^{-1}\Phi'.
\end{align*}
Note that $\Phi'$ satisfies equation (\ref{nadj}) and $\mathcal{S}P^{-1}\Phi'=0$. The latter follows from the condition $P(\theta)^{-1}\Phi'(\theta) ={\bf e}$ and ${\mathcal S}{\bf e}=\nu_1=0$. Finally, from the definition of $R(\theta)$, equation (\ref{Rtan}), we deduce that $\mathcal{S}P^{\top}(\theta)R(\theta)=0$ and hence
\[ P^T(R'+\omega_0^{-1}J^{\top}R)=0.\]
Since $P^T(\theta)$ is non-singular for all $\theta$, $R$ satisfies equation (\ref{adj}) and can thus be identified as the PRC.

\subsection{Variational principle}

Let $\beta(t)$ be the phase at time $t$. Since $\X(t)$ is constant between jumps, $\beta(t)$ is also taken to be constant between jumps. We can now state the variational principle for the stochastic phase \cite{Bressloff18a}: $\beta(t)$ for $t \in {\mathcal T}$ is determined by requiring $\beta(t)=b_t(\theta_t)$, where $b_t(\theta_t)$ for a prescribed time dependent weight $\theta_t$ is a local minimum of the following variational problem: 
\begin{equation}\label{minim}
\underset{b\in {\mathcal N}(b_t(\theta_t))} \inf\| \X(t)-\Phi(b)\|_{\theta_t} =\| \X(t)-\Phi(b_t(\theta_t))\|_{\theta_t} , \quad t \in {\mathcal T},
\end{equation}
with ${\mathcal N}\big(b_t(\theta_t)\big)$ denoting a sufficiently small neighborhood of $b_t\big(\theta_t\big)$. The minimization scheme is based on the orthogonal projection of the solution on to the limit cycle with respect to the weighted Euclidean norm at some $\theta_t$. Write the partial derivative of \eqref{minim} as
\begin{equation}\label{eq: Gi definition}
\G_0(\x,b,\theta):=\frac{\partial }{\partial b} \| \x-\Phi(b)\|^2_{\theta}  =-2\left \langle \x-\Phi(b),\Phi'(b)\right \rangle_{\theta}.
\end{equation}
At the minimum, 
\begin{equation}\label{eq: minimum G 0}
\G_0\big(\X(t),\beta(t),\theta_t\big)=0.
\end{equation}
We stipulate that the location of the weight must coincide with the location of the minimum, i.e. $\beta(t) = \theta_t$, so that $\beta(t)$ must satisfy the implicit equation
\begin{equation}\label{eq: minimum G 2}
\G\big(\X(t),\beta(t)\big):=\G_0\big( \X(t) ,\beta(t),\beta(t)\big)=0.
\end{equation}
Let $\tau_k$ denote the time of the $k$th reaction in a given realization of the stochastic CRN. (For the moment, we don't specify which reaction occurs.) Since $\X(t)$ is constant for $t\in [\tau_k,\tau_{k+1})$, $\beta(t)$ is also constant over this interval. It is clear that there could be more than one minimum of \eqref{minim}, and therefore this equation may be insufficient to determine $\beta(\tau_{k+1})$. We therefore further stipulate that $\beta(\tau_{k+1})$ is the local minimum of \eqref{eq: Gi definition} that also minimizes $\big| \beta(\tau_{k+1}) -\beta(\tau_{k+1}^-)\big|$. In other words, we emphasize continuity over finding the global minimum. This is consistent with our  amplitude/phase decomposition, as well as the classical isochron phase decomposition of stochastic oscillators. There are several reasons for this choice. We wish to identify a linear phase equation, which requires the increments in $\beta(t)$ to be small. Furthermore, situations where a global minimum switches to being a local minimum, and the corresponding jump in phase is much greater than $O(\Omega^{-1})$, are relatively rare. Trying to understand this highly nonlinear phenomenon would take us beyond the scope of this paper.

\subsection{Linear phase approximation}

As it stands, the variational principle generates an exact, implicit equation for the stochastic phase. Here we carry out a perturbation expansion to show that the leading order dynamics of the phase over small time intervals recovers the explicit isochronal phase equation (\ref{theta}); we will refer to this as the {\em linear phase approximation}.
 Let $\delta \X_k = \X(\tau_{k+1}) - \X(\tau_k)$. It is immediate from our earlier equations that $\big| \delta \X_k\big| \leq \Omega^{-1}\sup_{1\leq a \leq M}\norm{S_a}$. We use \eqref{eq: minimum G 2} to estimate $\delta \beta_k$ to leading order. It must be that
\begin{equation}
\G_0\big( \X(\tau_{k+1}) ,\beta(\tau_{k+1}),\beta(\tau_{k+1})\big)-\G_0\big( \X(\tau_k) ,\beta(\tau_k),\beta(\tau_k)\big) = 0.
\end{equation}
This means that, to leading order,
\begin{align}
\frac{\partial}{\partial \X}\G_0\big( \X(\tau_k) ,\beta(\tau_k),\beta(\tau_k)\big)\delta \X_k  + \frac{\partial}{\partial \beta}\G_0\big( \X(\tau_k) ,\beta(\tau_k),\beta(\tau_k)\big)\delta \beta_k  = 0.
\end{align}
Rearranging, we find that
\begin{equation}\label{eq: delta beta k}
\delta \beta_k \simeq \frac{1}{\mathfrak{M}(\X_k,\beta_k)} \bigg\langle \Phi'\big(\beta(\tau_k)\big),\delta \X_k \bigg\rangle_{\beta(\tau_k)},
\end{equation}
where
\begin{align}
\mathfrak{M}(\X,\beta) :=\frac{\partial \G(\X,\beta)}{\partial \beta}=- \bigg\langle \Phi''(\beta),\X - \Phi(\beta) \bigg\rangle_{\beta(\tau_k)}+ \bigg\langle \Phi'(\beta), \Phi'(\beta) \bigg\rangle_{\beta(\tau_k)} \\-  \bigg\langle \Phi'(\beta),\frac{d}{d\beta}\big(P(\beta)P(\beta)\big)^{-T}\big(\X - \Phi(\beta)\big) \bigg\rangle.\nonumber
\end{align}
\james{It follows from our earlier definitions that $\big\langle \Phi'(\beta), \Phi'(\beta) \big\rangle_{\beta(\tau_k)} = 1$. Now as long as $\norm{\X - \Phi(\beta)} \ll 1$ (we determine the probability of this in $\S$ 5),  $\mathfrak{M}(\X,\beta) \geq \frac{1}{2}$ (say). We can then use Taylor's Theorem to infer that, since $\delta \X_k = O(\Omega^{-1})$, the error in \eqref{eq: delta beta k} is $O(\Omega^{-2})$.} Note that, since by assumption $\beta$ is a local minimum, $\mathfrak{M}(\X,\beta)$ must be positive, because it is the curvature at the local minimum.

If we now take note of the reaction label $a$, it follows by analogy with \eqref{rtc} that
\begin{align}
\label{ba}
\beta(t) = \beta(0)+\frac{1}{\Omega}\sum_{a=1}^M \int_0^t \frac{1}{\mathfrak{M}(\X(s),\beta(s))} \bigg\langle \Phi'\big(\beta(s)\big), {\bf S}_a \bigg\rangle_{\beta(s)} d {Y}^a(s)+ O\bigg(  \Omega^{-2}{\mathcal M}_{0,t}\bigg).
\end{align}
Here 
\[
 {Y}^a(s) = Y^a \bigg(\Omega\int_0^s \lambda^a\big(\X(u)\big)du \bigg),
\]
which reduces to equation (\ref{Ya}) on taking $\X(u)=\Phi(\beta(u))$.
We note that the above Poisson integral is, by definition, just the summation over the set of jumps, see section 3. That is, by definition,
\begin{align*}
&\int_0^t \frac{1}{\mathfrak{M}(\X(s),\beta(s))} \bigg\langle \Phi'\big(\beta(s)\big),  {\bf S}_a \bigg\rangle_{\beta(s)} d\bar{Y}^a(s)  = \Omega^{-1}\sum_{k>0} \frac{H(t-t_k^a) }{\mathfrak{M}(\X(t^a_k),\beta(t^a_k))} \bigg\langle \Phi'\big(\beta(t^a_k)\big), {\bf S}_a \bigg\rangle_{\beta(t^a_k)}, 
\end{align*}
where $\lbrace t^a_i \rbrace_{i\in \mathbb{Z}^+}$ are the jump times of $Y^a \bigg(\Omega \int_0^s \lambda^a\big(\X(u)\big)du \bigg)$.

We wish to represent the above equation as a time-integral plus a \textit{martingale}. (A martingale is a stochastic process whose increments over small time-intervals are independent of the past history, and of mean zero \cite{Protter05,Anderson15}. For example any stochastic integral with respect to a Wiener Process is a martingale. \james{The martingales in the following equations are often referred to as \textit{compensated Poisson processes}.}) In this way, the form of the equation will resemble that of a standard stochastic differential equation. Another important reason that we wish to put the equation in this form is that it will allow us to use a martingale inequality to obtain bounds on the system leaving a neighborhood of the limit cycle, see \S 5. To this end, equation (\ref{ba}) can be rewritten as
\begin{align}
\label{better}
\beta(t) = \beta(0)+\sum_{a=1}^M \int_0^t \frac{1}{\mathfrak{M}(\X(s),\beta(s))} \bigg\langle \Phi'\big(\beta(s)\big),  {\bf S}_a \bigg\rangle_{\beta(s)}\lambda_a\big( \X(s) \big)ds   + \Omega^{-1}\widetilde{Y}(t) +O\bigg( \Omega^{-2}{\mathcal M}_{0,t}\bigg),
\end{align}
where $\widetilde{Y}(t)$ is the martingale
\begin{align}\label{phase martingale}
\widetilde{Y}(t) =  \sum_{a=1}^M \int_0^t \frac{1}{\mathfrak{M}(\X(s),\beta(s))} \bigg\langle \Phi'\big(\beta(s)\big),  {\bf S}_a \bigg\rangle_{\beta(s)} \big\lbrace d{Y}^a(s) - \Omega \lambda_a\big( \X(s) \big)ds \rbrace .
\end{align}
(Note that $\widetilde{Y}(t)$ is discontinuous, unlike the stochastic integral.)

Equation (\ref{better}) has the same basic structure as the isochronal phase equation (\ref{theta}). Indeed, it reduces to equation (\ref{theta}) if we set $\X(s)=\Phi(\beta(s))$ on the right-hand side. This follows from the observations
\[\M(\Phi(\beta),\beta)=\bigg\langle \Phi'(\beta), \Phi'(\beta) \bigg\rangle_{\beta}=\bigg\langle P^{-1}(\beta) \Phi'(\beta),P^{-1}(\beta) \Phi'(\beta) \bigg\rangle=1,\]
and
\[\bigg\langle \Phi'\big(\beta\big),  {\bf S}_a \bigg\rangle_{\beta(s)}=\bigg\langle [P(\theta)P^{\top}(\theta)]^{-1}\Phi'(\theta), {\bf S}_a\bigg\rangle=\sum_{j=1}^K S_{ja}R_j(\theta),\]
on using equation (\ref{Rtan}). In particular,
\begin{align*}
\sum_{a=1}^M\frac{1}{\mathfrak{M}(\X(s),\beta(s))} \bigg\langle \Phi'\big(\beta(s)\big),  {\bf S}_a \bigg\rangle_{\beta(s)}\lambda_a\big( \X(s) \big) = \omega_0 + O\big( \norm{ \X(s) - \Phi(\beta(s))}^2 \big). 
\end{align*}
We thus obtain the leading order phase equation \james{
\begin{align}\label{eq: linear beta}
\beta(t) \simeq \beta(0) + \omega_0 t+ \Omega^{-1}\widetilde{Y}(t) +  O\bigg( \mathcal{M}_{0,t} \Omega^{-2}+ t\sup_{s\in [0,t]}\norm{ \X(s) - \Phi(\beta(s))}^2 \bigg).
\end{align}}
Furthermore, this linearized phase is typically accurate over timescales of $O(\Omega)$, \james{ as long as $t\sup_{s\in [0,t]}\norm{ \X(s) - \Phi(\beta_s)}^2$} remains small (we are going to bound this in the next section). This is because, as we saw in section 3, the total number of reactions over the time interval $[0,t]$ is of order $\Omega\overline{\lambda}t$, with very high probability.  For large $\Omega$, the probability law of the martingale $\widetilde{Y}(t)$ may be approximated by the stochastic integral
\begin{align*}
\widetilde{Y}(t) \simeq  \sqrt{\Omega}\sum_{a=1}^M  \int_0^t \frac{1}{\mathfrak{M}(\X(s),\beta(s))} \bigg\langle \Phi'\big(\beta(s)\big),  {\bf S}_a \sqrt{\lambda_a(\X(s))} \bigg\rangle_{\beta(s)}dW^a(s),
\end{align*}
where $\lbrace W^a(s) \rbrace$ are independent Wiener processes \cite{Kurtz71,Kurtz76}. In this way, we see that the probability law of the phase equation (\ref{eq: linear beta}) may be approximated to be
\begin{align*}
\beta(t) \simeq \beta(0) + t\omega_0 + \frac{1}{\sqrt{\Omega}}\sum_{a=1}^M  \int_0^t \frac{1}{\mathfrak{M}(\Phi(\beta(s)),\beta(s))} \bigg\langle \Phi'\big(\beta(s)\big), {\bf S}_a \sqrt{\lambda_a(\Phi(\beta(s)))} \bigg\rangle_{\beta(s)}dW^a(s).
\end{align*}
Note that this agrees with the diffusion approximation for the phase in equation  \eqref{phaseFP}.

\subsection{Example: the Brusselator} One benchmark model of a biochemical oscillator is the Brusselator, which is an idealized model of an autocatalytic reaction, in which at least one of the reactants is also a product of the reaction \cite{Epstein98}. The model consists of two chemical species $X$ and $Y$ interacting through the following reaction scheme:
\begin{align*}
\emptyset &\overset{a} \rightarrow X\\
X &\overset{b} \rightarrow Y\\
2X+Y &\overset{c} \rightarrow 3X\\
X &\overset{d} \rightarrow \emptyset
\end{align*}
These reactions describe the production and degradation of an $X$ molecule, an $X$ molecule spontaneously transforming into a $Y$ molecule, and two molecules of $X$ reacting with a single molecule of $Y$ to produce three molecules of $X$. The corresponding mass-action kinetic equations for $u_1=[X],u_2=[Y]$ are (after rescaling so that $c=d=1$)
\begin{subequations}
\label{Brussel}
\begin{align}
\frac{du_1}{dt}&=a-(b+1)u_1+u_1^2u_2,\\
\frac{du_2}{dt}&=bu_1-u_1^2u_2.
\end{align}
\end{subequations}
The system has a fixed point at $u_1^*=a,u_2^*=b/a$, which is stable when $b<a^2+1$ and unstable when $b>a^2+1$. Moreover, the fixed point undergoes a Hopf bifurcation at the critical value $b=a^2+1$ for fixed $a$, leading to the formation of a stable limit cycle.

We simulated the stochastic Brusselator with $\Omega = 3000$, $a=c=d=1$ and $b=2.5$. The results are plotted in Figure \ref{Fig:Brusselator}. It can be seen that the phase under the linear phase approximation (linear phase) is initially in close agreement with the exact variational phase. However, after one orbit of the limit cycle, there is already a significant difference between the two. Despite this nonlinear diffusion of the phase, the system still stays in close proximity to the deterministic limit cycle for the entire length of the simulation. This can be explained using the analysis of \S 5, where we demonstrate that with very high probability the system stays close to the limit cycle for times of $O\big(\exp(b\Omega C)\big)$, where $b$ is the rate of decay towards the limit cycle, and $C$ is a constant.

\begin{figure}[t!]\label{Fig:Brusselator}
\begin{center}
\includegraphics[width = 14cm]{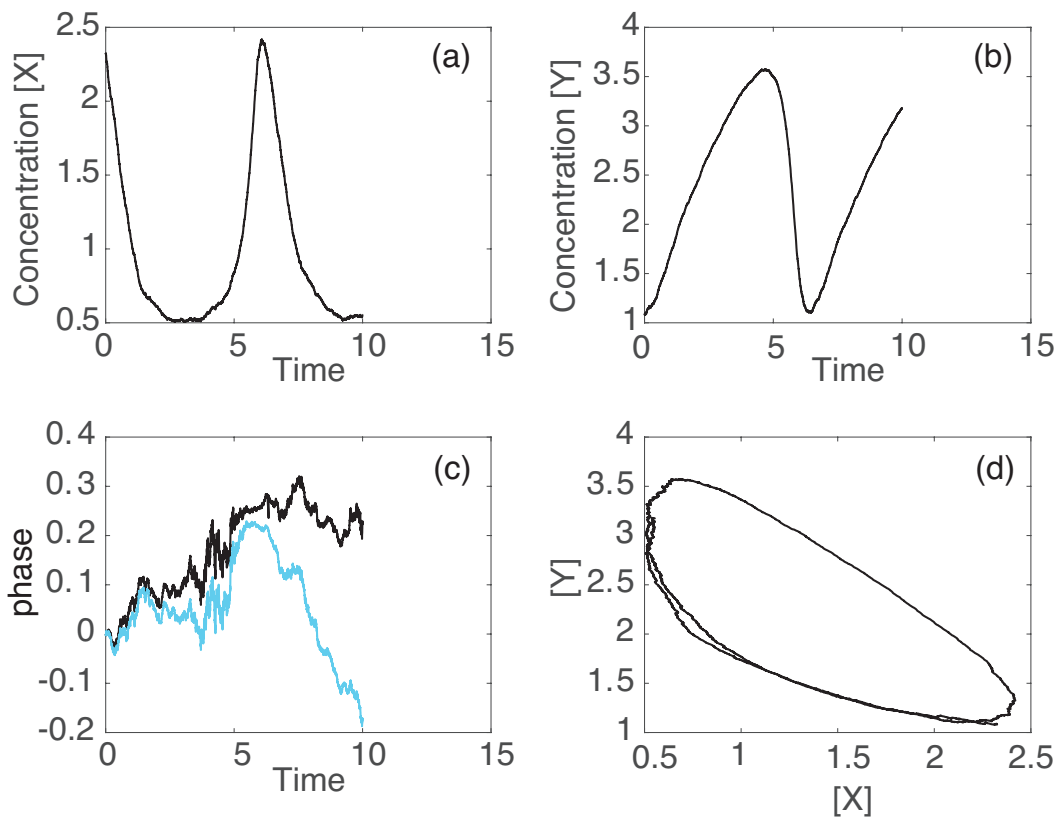}
\caption{Numerical simulation of the stochastic Brusselator, with $a=c=d=1$, $b=2.5$ and the total number of particles $\Omega=3000$. Initially, $[X]$ and $[Y]$ are taken to lie exactly on the limit cycle. (a) Plot of concentration $[X]$ as a function of time $t$. (b) Corresponding plot for $[Y]$. (c) Plot of the linear phase (light trajectory) and the exact variational phase (dark trajectory). (We have subtracted $\omega_0 t$ off both solutions.) (d) Phase portrait in $[X]-[Y]$ plane. It is notable that even with such a large $\Omega$, the linear phase is substantially different from the variational phase after just one orbit of the limit cycle.}
\label{fig1}
\end{center}
\end{figure}

\section{Bounding the probability of the CRN leaving a neighborhood of the limit cycle}

The probability bounds of this section are accurate in the limit that $\Omega b \gg 1$ (recall that $b$ is the rate of decay towards the limit cycle in the deterministic dynamics). Throughout this analysis, we assume that
\begin{align}
\sup_{\theta\in [0,2\pi]} \max\bigg\lbrace \norm{\mathbf{P}(\theta)}, \norm{\bf{P}(\theta)^{-1}}\bigg\rbrace = O(1) \label{eq: P assumption 1} \\
\sup_{\bf{Z} \in \mathbb{R}^K, 1\leq i,j,k \leq K}\bigg| \frac{\partial^2}{\partial Z^j \partial Z^k}F_i({\bf Z})\bigg| = O(1).
\end{align}
\james{It should also be noted that $\omega_0$ can be arbitrarily large and the bounds would still be accurate (this is precisely what was observed in the numerical simulation of the chemical Brusselator earlier - the diffusion of the phase was already substantial after one orbit of the limit cycle, however the system still stayed close to the limit cycle).}
Define the amplitude ${\bf v}$ and weighted amplitude ${\bf w}$ according to
\begin{align}
{\bf v}(t) =  \X(t) - \Phi\big( \beta(t) \big) ,\quad 
{\bf w}(t)= P\big(\beta(t) \big)^{-1}{\bf v}(t) .
\end{align}
As is explained in more detail in \cite{Bressloff18a}, the reason that we make this change of variable is that, to leading order, $\norm{\bf{w}(t)}$ decays uniformly in time.
The main result that we prove in this section is that there exists a constant $C > 0$ such that
\begin{equation}\label{eq: main result}
\mathbb{P}\bigg(\sup_{t\in [0,T]}\norm{{\bf w}(t)} \geq \zeta \bigg) \leq  bT\bigg\lbrace \exp\big(-\Omega b\zeta^2 C\big) \bigg\rbrace,
\end{equation}
for all sufficiently small $\zeta$. We assume throughout this section that
\begin{equation}
\norm{{\bf w}(0)} \leq \frac{\zeta}{2}.
\end{equation}
As in \S 4.3, suppose that $\tau_k$ is the time of the $k$th reaction. To leading order in $\delta \beta_k , \delta X_k$, \james{
\begin{align}
{\bf w}(\tau_{k+1}) - {\bf w}(\tau_k) \simeq \frac{\partial}{\partial \beta}\bigg\lbrace P\big(\beta(\tau_k) \big)^{-1}\bigg(\X(\tau_k) - \Phi\big( \beta(\tau_k) \big) \bigg) \bigg\rbrace \delta \beta_k + P\big(\beta(\tau_k) \big)^{-1}\delta \X_k,
\end{align}
and the error in the above approximation is $O\big( \norm{\delta \X_k}^2, \delta \beta_k^2 \big) = O\big( \Omega^{-2}\big)$.}

Note that
\begin{align}
 \frac{\partial}{\partial \beta}\bigg\lbrace P\big(\beta \big)^{-1}\bigg(\X(t) - \Phi\big(\beta(t)  \big) \bigg) \bigg\rbrace
 &= - \omega_0^{-1} P\big(\beta \big)^{-1}J(\beta)\bigg(\X(t) - \Phi\big( \beta \big) \bigg) +\omega_0^{-1} \mathcal{S}{\bf w}(t)  \nonumber \\&\quad  -P\big(\beta \big)^{-1}\Phi'\big(\beta\big),
 \end{align}
which follows from the identity 
\[
\omega_0\frac{d}{dt}P\big(t\big) = J(t) P(t) - P(t) \mathcal{S}(t).
\]
Substituting \eqref{rtc} and \eqref{eq: linear beta}, we thus find that for any $u\in [0,T]$ and $t > u$,
\begin{align}
{\bf w}(t) &= {\bf w}(u)+  O\bigg(\Omega^{-2}\mathcal{M}_{u,t} + (t-u)\norm{\bf{v}(t)}^2 \bigg)
+ \Omega^{-1}[ \breve{Y}(t) - \breve{Y}(u)] \\
&\quad + \int_u^t\bigg\lbrace \mathcal{S}{\bf w}(s) + P\big(\beta(s)\big)^{-1}\bigg( F\big( \X(s) \big) - F\big(\Phi(\beta(s))\big)-  J\big(\beta(s)\big)\big\lbrace \X(s) - \Phi\big( \beta(s) \big) \big\rbrace   \bigg) \bigg\rbrace ds  .\nonumber
\end{align}
where $\mathcal{M}_{u,t}$ is, as before, the total number of reactions over the time interval $[u,t]$. We have used the fact that $\omega_0 \Phi'(\beta)=F\big(\Phi(\beta)\big)$ and the increments in $\lbrace \X,\beta,{\bf w} \rbrace$ are all $O(\Omega^{-1})$ (which we noted in the previous section). Here $\breve{Y}(t)$ is the Martingale
\begin{align}
&\breve{Y}(t) = \sum_{a=1}^M \int_0^t P\big(\beta(s)\big)^{-1} {\bf S}_a d\grave{Y}_a(s)\\
& \quad +\int_0^t \bigg\lbrace- \omega_0^{-1} P\big(\beta(s) \big)^{-1}J(\beta(s))\bigg(\X(s) - \Phi\big( \beta(s) \big) \bigg) +\omega_0^{-1} \mathcal{S}{\bf w}(s)  -P\big(\beta(s) \big)^{-1}\Phi'\big(\beta(s)\big)\bigg\rbrace d\widetilde{Y}(s) \nonumber 
\end{align}
By Taylor's Theorem, 
\begin{align}
F\big( \X(s) \big) -  J\big(\beta(s)\big)\big\lbrace \X(s) - \Phi\big( \beta(s) \big) \big\rbrace - F\big(\Phi(\beta(s))\big)
= O\bigg( \norm{\X(s) - \Phi(\beta_s)}^2 \bigg).
\end{align}
Since $\norm{{\bf v}(t)} = O\big(\norm{{\bf w}(t)}\big)$, noting the assumption in \eqref{eq: P assumption 1}, this means that
\begin{align}
{\bf w}(t) &= {\bf w}(u)+ \int_u^t \mathcal{S}{\bf w}(s)+  O\bigg(\Omega^{-2}\mathcal{M}_{u,t} + (t-u)\sup_{r\in [u,t]}\norm{{\bf w}(r)}^2 \bigg)
+ \Omega^{-1}[ \breve{Y}(t) - \breve{Y}(u)] \nonumber
\end{align}
Finally, we apply the changes of variable ${\bf w}(t) \to \big\langle {\bf w}(t),{\bf w}(t)\big\rangle \to \sqrt{\big\langle {\bf w}(t),{\bf w}(t)\big\rangle}= \norm{{\bf w}(t)}$, and find that
\begin{align}
&\norm{{\bf w}(t)} = \norm{{\bf w}(u)} + \int_u^t \frac{1}{\norm{{\bf w}(s)}} \bigg\langle {\bf w}(s),\mathcal{S}{\bf w}(s) \bigg\rangle ds+O\bigg(\Omega^{-2}\mathcal{M}_{u,t} + (t-u)\sup_{r\in [u,t]}\norm{{\bf w}(r)}^2 \bigg)\nonumber \\
&\qquad 
+\frac{1}{\Omega} \int_u^t\frac{1}{\norm{{\bf w}(s)}}\big\langle {\bf w}(s), d\breve{Y}(s)\big\rangle.
\end{align}

Furthermore our stability assumption on the limit cycle, see equation (\eqref{eq: nu bound}), implies that
\[
\big\langle {\bf w}(s),\mathcal{S}{\bf w}(s)\big\rangle \leq -b\norm{{\bf w}(s)}^2.
\]
Note that the first element of ${\bf w}$ is identically zero, as a direct consequence of its definition, which is what leads to the strict decay in the above equation.

Hence, as long as $\norm{w(s)} \leq \frac{b}{2}$,
\begin{align}
\norm{{\bf w}(t)} \leq \norm{{\bf w}(u)} -\frac{b}{2} \int_u^t \norm{{\bf w}(s)} ds
+O\big(\mathcal{M}_{u,t}\Omega^{-2} \big) + \mathcal{Y}(t) - \mathcal{Y}(u),\label{eq: norm w inequality}
\end{align}
where
\[
\mathcal{Y}(t) = \frac{1}{\Omega}\int_0^t\frac{1}{\norm{{\bf w}(s)}}\big\langle {\bf w}(s), d\breve{Y}(s)\big\rangle.
\]
We now outline a set of events that, taken together, ensure that $\sup_{t \leq T}\norm{{\bf w}(t)} \leq \zeta$. This means that, to bound the probability that 
$\sup_{t \leq T}\norm{{\bf w}(t)} > \zeta$, it suffices to bound each of the following events. These events are
\begin{subequations}
\begin{align}
\norm{{\bf w}(0)} \leq \frac{\zeta}{2}\label{eq: event 1} \\
\sup_{t\in [0,T]} \sup_{s \in [0,b^{-1}]}\bigg\lbrace  \mathcal{Y}(t+s) - \mathcal{Y}(t) \bigg\rbrace \leq \frac{\zeta}{8} \\
\sup_{s\in [0,T]} \sup_{t \in [0,b^{-1}]}\Omega^{-2}\mathcal{M}_{s,s+t} \leq \frac{\zeta}{8}.\label{eq: event 3} 
\end{align}
\end{subequations}
The timescale of $b^{-1}$ has been chosen for a reason: it is the timescale of the decay towards the limit cycle. To see why the above three events are sufficient to ensure that $\sup_{t \leq T}\norm{{\bf w}(t)} \leq \zeta$, let us suppose for a contradiction that \eqref{eq: event 1}-\eqref{eq: event 3} all hold, but that
$\norm{{\bf w}(t)} = \zeta$ for some $t\in [0,T]$, and that $\norm{{\bf w}(s)} < \zeta$ for all $s\in [0,t)$ \james{(i.e. we are taking $t$ to be the first time that $\norm{{\bf w}}$ attains $\zeta$)}. Let $\tau  = \sup\big\lbrace u\leq t \; : \; \norm{{\bf w}(u)} = \frac{\zeta}{2} \big\rbrace$, noting that $\tau$ exists since, by assumption, $\norm{{\bf w}(0)}\leq {\zeta}/{2}$. 

Suppose for a contradiction that $t-\tau \geq b^{-1}$. Then, since $\norm{{\bf w}(s)}\geq {\zeta}/{2}$ for all $s\in [\tau,t]$, it follows from \eqref{eq: norm w inequality} that
\begin{align*}
\norm{{\bf w}(t)} &\leq \norm{{\bf w}(t-b^{-1})} -\frac{b}{2b}\frac{\zeta}{2}
+\frac{\zeta}{8} + \frac{\zeta}{8} 
\leq \norm{{\bf w}(t-b^{-1})}.
\end{align*}
This contradicts our assumption that $\norm{{\bf w}(s)} < \zeta$ for all $s \in [0,t)$. We can thus conclude that $t-\tau < b^{-1}$. But in this case, it is immediate that 
\begin{align*}
\norm{{\bf w}(t)} \leq & \norm{{\bf w}(\tau)} + \frac{\zeta}{8} + \frac{\zeta}{8}
\leq  \frac{3\zeta}{4},
\end{align*}
which contradicts our assumption that $\norm{{\bf w}(t)} = \zeta$. 

We can therefore conclude that the events in \eqref{eq: event 1} -\eqref{eq: event 3} imply that $\sup_{s\in [0,T]}\norm{{\bf w}(s)} < \zeta$. Define the event 
\begin{equation}
\mathcal{H}_{\zeta} = \bigg\lbrace \sup_{s\in [0,T]}\norm{{\bf w}(s)} \leq \zeta \bigg\rbrace.
\end{equation}
Given that, by stipulation, $\norm{{\bf w}(0)} \leq {\zeta}/{2}$, this means that
\begin{align}
 \label{eq: to prove P decomposition}
&\mathbb{P}\bigg(\sup_{s\in [0,T]}\norm{{\bf w}(s)} \geq \zeta \bigg) \\ &\quad \leq \mathbb{P}\bigg(\sup_{t\in [0,T]} \sup_{s \in [0,b^{-1}]}\bigg\lbrace  \mathcal{Y}(t+s) - \mathcal{Y}(t) \bigg\rbrace > \frac{\zeta}{8} \text{ and }\mathcal{H}_{\zeta} \bigg)
+ \mathbb{P}\bigg(\sup_{s\in [0,T]} \sup_{t \in [0,b^{-1}]}\Omega^{-2}\mathcal{M}_{s,s+t} > \frac{\zeta}{8} \bigg).\nonumber
\end{align} 
Note that we have appended the event $\mathcal{H}_a$ into the first term on the right, because, as has just been demonstrated, this must always hold as long as the other events hold as well. We now bound the above two probabilities, showing that they are each $O\big(\exp(b\Omega C)\big)$ for some constant $C > 0$.

Let $t_i = {i}/{2b}$. First, using a union of events bound,
\begin{equation}
\mathbb{P}\bigg(\sup_{s\in [0,T]} \sup_{t \in [0,b^{-1}]}\Omega^{-2}\mathcal{M}_{s,s+t} > \frac{\zeta}{8} \bigg) 
\leq  \sum_{i=0}^{\lfloor 2bT \rfloor} \mathbb{P}\bigg(\sup_{t \in [t_i,t_{i+1}]}\Omega^{-2}\mathcal{M}_{s,s+t} > \frac{\zeta}{16} \bigg).
\end{equation}
Furthermore, it follows from \eqref{eq: bound number of reactions} that for some constant $C$,
\[
\mathbb{P}\bigg( \sup_{t \in [t_i,t_{i+1}]}\Omega^{-2}\mathcal{M}_{s,s+t} > \frac{\zeta}{8} \bigg) \leq \exp\big(-\Omega bC\big),
\]
as long as $b\Omega$ is sufficiently large. We thus find that
\[
 \mathbb{P}\bigg(\sup_{s\in [0,T]} \sup_{t \in [0,b^{-1}]}\Omega^{-2}\mathcal{M}_{s,s+t} > \frac{\zeta}{8} \bigg) = O\bigg(Tb^{-1}\exp\big(Cb\Omega\big) \bigg),
\]
as required.

We turn to the second inequality in \eqref{eq: to prove P decomposition}. We first claim that
\begin{align}
& \mathbb{P}\bigg(\sup_{t\in [0,T]} \sup_{s \in [0,b^{-1}]}\bigg\lbrace  \mathcal{Y}(t+s) - \mathcal{Y}(t) \bigg\rbrace > \frac{\zeta}{8} \text{ and }\mathcal{H}_{\zeta}\bigg)\nonumber \\
&\qquad  \leq \sum_{i=0}^{\lfloor 2bT \rfloor}  \mathbb{P}\bigg( \sup_{s \in [t_i,t_{i+1}]}\bigg|  \mathcal{Y}(s) - \mathcal{Y}(t_i) \bigg|  > \frac{\zeta}{24} \text{ and }\mathcal{H}_{\zeta} \bigg).\label{eq: to prove P decomposition 2}
\end{align}
To see why this is true, suppose that $t+s\in [t_i,t_{i+1}]$, and $t\in [t_j,t_{j+1}]$, so that either $j=i$ or $j=i-1$. Then
\begin{align*}
\bigg|  \mathcal{Y}(t+s) - \mathcal{Y}(t) \bigg|  \leq \bigg|  \mathcal{Y}(t+s) - \mathcal{Y}(t_{i}) \bigg|
+\bigg|  \mathcal{Y}(t_i) - \mathcal{Y}(t_j) \bigg| + \bigg|  \mathcal{Y}(t_j) - \mathcal{Y}(t) \bigg|.
\end{align*}
Hence if each of the terms on the right is less than ${\zeta}/{24}$, their sum must be less than ${\zeta}/{8}$. This establishes \eqref{eq: to prove P decomposition 2}. It thus suffices for us to bound
\[
\mathbb{P}\bigg( \sup_{s \in [t_i,t_{i+1}]}\bigg|  \mathcal{Y}(s) - \mathcal{Y}(t_i) \bigg|  > \frac{\zeta}{24}\text{ and }\mathcal{H}_{\zeta} \bigg).
\]
Because the exponential function is positive and increasing,
\begin{align*}
&\mathbb{P}\bigg( \sup_{s \in [t_i,t_{i+1}]}\bigg|  \mathcal{Y}(s) - \mathcal{Y}(t_i) \bigg|  > \frac{\zeta}{24} \text{ and }\mathcal{H}_{\zeta} \bigg) \\
&\leq \mathbb{P}\bigg( \sup_{s \in [t_i,t_{i+1}]} \exp\bigg(\alpha \bigg\lbrace  \mathcal{Y}(s) - \mathcal{Y}(t_i)\bigg\rbrace\bigg) > \exp\bigg\lbrace\frac{\zeta \alpha}{24}\bigg\rbrace \text{ and }\mathcal{H}_{\zeta} \bigg) \\
&+ \mathbb{P}\bigg( \sup_{s \in [t_i,t_{i+1}]} \exp\bigg(-\alpha \bigg\lbrace  \mathcal{Y}(s) - \mathcal{Y}(t_i)\bigg\rbrace\bigg) > \exp\bigg\lbrace\frac{\zeta\alpha}{24}\bigg\rbrace\text{ and }\mathcal{H}_{\zeta} \bigg).
\end{align*}
Since $\mathcal{Y}$ is a martingale, both $\exp\big(\alpha \big\lbrace \mathcal{Y}(t) -\mathcal{Y}(t_i)\big\rbrace \big)$ and $\exp\big(-\alpha \big\lbrace \mathcal{Y}(t) -\mathcal{Y}(t_i)\big\rbrace \big)$ are submartingales for any positive constant $\alpha$, \james{since the exponential function is convex (and we have also used Jensen's Inequality). We can thus use Doob's submartingale inequality \cite[A.2.1]{Anderson15}}  to conclude that
\begin{align}
\mathbb{P}\bigg( \sup_{s \in [t_i,t_{i+1}]} \exp\bigg(\alpha \bigg\lbrace  \mathcal{Y}(s) - \mathcal{Y}(t_i)\bigg\rbrace\bigg) > \exp\bigg(\frac{\zeta\alpha}{24}\bigg)\text{ and }\mathcal{H}_{\zeta} \bigg)\nonumber  \\
\leq \mathbb{E}\bigg[ \exp\bigg(\alpha \bigg\lbrace  \mathcal{Y}(t_{i+1}) - \mathcal{Y}(t_i) - \frac{\zeta}{24} \bigg\rbrace\bigg)\mathbf{1}\big\lbrace \mathcal{H}_{\zeta} \big\rbrace \bigg],\label{eq final inequality 1}
\end{align}
and
\begin{align}
\mathbb{P}\bigg( \sup_{s \in [t_i,t_{i+1}]} \exp\bigg(-\alpha \bigg\lbrace  \mathcal{Y}(s) - \mathcal{Y}(t_i)\bigg\rbrace\bigg) > \exp\bigg(\frac{\zeta\alpha}{24}\bigg)\text{ and }\mathcal{H}_{\zeta} \bigg)\nonumber \\
\leq \mathbb{E}\bigg[ \exp\bigg(-\alpha \bigg\lbrace  \mathcal{Y}(t_{i+1}) - \mathcal{Y}(t_i) + \frac{\zeta}{24} \bigg\rbrace\bigg)\mathbf{1}\big\lbrace \mathcal{H}_{\zeta} \big\rbrace \bigg].\label{eq: to bound final 2}
\end{align}
We now bound the first of the above expectations. Retracing the above derivation, we can write
\begin{equation}\label{eq: mathcal Y t}
\mathcal{Y}_t -\mathcal{Y}_s= \Omega^{-1}\sum_{a=1}^M \int_s^t g_a\big(r, {\bf X}(r)\big) d\mathcal{N}_a(r) - \int_s^t g_a\big( r, {\bf X}(r)\big)\lambda_a(r) dr,
\end{equation}
for functions $\lbrace g_a \rbrace_{a=1}^M$ that are uniformly bounded in the tube $\norm{ {\bf X}(t) - \Phi(t)}_{\beta(t)} \leq \zeta$. Let this uniform bound be $\bar{g}$. Define
\begin{equation}
\mathcal{H}_{\zeta,t} = \bigg\lbrace \sup_{s\in [0,t]}\norm{{\bf w}(s)} \leq \zeta \bigg\rbrace,
\end{equation}
and note that $\mathcal{H}_{\zeta} \subseteq \mathcal{H}_{\zeta ,t}$. Using \eqref{eq: mathcal Y t} and standard properties of the Poisson process \cite{Anderson15},
\begin{multline*}
\frac{d}{dt}\mathbb{E}\bigg[\exp\big(\alpha \mathcal{Y}_t -\alpha \mathcal{Y}_s \big)\mathbf{1}\big\lbrace \mathcal{H}_{\zeta,t} \big\rbrace \bigg]  \\ \leq \Omega \mathbb{E}\bigg[\exp\big(\alpha \mathcal{Y}_t-\alpha \mathcal{Y}_s\big)   \sum_{a=1}^M \lambda_a(t)  \bigg\lbrace \exp\big( \alpha \Omega^{-1} g_a\big(t,{\bf X}(t)\big) \big) - 1 - \alpha \Omega^{-1} g_a\big( t,{\bf X}(t)\big)\bigg\rbrace \mathbf{1}\big\lbrace \mathcal{H}_{\zeta,t} \big\rbrace\bigg]  \\
\leq \Omega \mathbb{E}\bigg[\exp\big(\alpha \mathcal{Y}_t\big)\mathbf{1}\big\lbrace \mathcal{H}_{\zeta,t} \big\rbrace\bigg]M e\bar{\lambda} \big(\alpha\Omega^{-1}\bar{g}\big)^2,
\end{multline*}
assuming that $\alpha$ is sufficiently small that $\alpha \Omega \bar{g} \leq 1$, so that we can use the second order Taylor bound for the exponential. By Gronwall's Inequality, we find that
\begin{equation}
\mathbb{E}\bigg[\exp\big(\alpha \mathcal{Y}_t  - \alpha\mathcal{Y}_s\big) \mathbf{1}\big\lbrace \mathcal{H}_{\zeta,t} \big\rbrace\bigg] \leq \exp\bigg( (t-s)eM \bar{\lambda}\big(\alpha\Omega^{-1}\bar{g}\big)^2\Omega \bigg).
\end{equation}
To optimize the bound in \eqref{eq final inequality 1}, and noting that $t-s = b^{-1}$, we thus choose
\[
\alpha = \Omega b \bigg(2eM \bar{\lambda}\bar{g}^2\bigg)^{-1} \frac{\zeta}{24} ,
\]
and we obtain that
\begin{align}
\mathbb{P}\bigg( \sup_{s \in [t_i,t_{i+1}]} \bigg\lbrace  \mathcal{Y}(s) - \mathcal{Y}(t_i)\bigg\rbrace >\frac{\zeta\alpha}{24} \text{ and }\mathcal{H}_{\zeta} \bigg)
\leq \exp\bigg( - C\Omega b \zeta^2  \bigg),
\end{align}
for some positive constant $C$. The bound for \eqref{eq: to bound final 2} is analogous.

\section{Discussion and Conclusion}
In this paper we outlined a variational definition of the phase of a stochastic chemical reaction oscillator. This definition is versatile and adaptable to a number of applications. First, we demonstrated that under a linear phase approximation, we recover the phase predicted by the diffusion approximation of the underlying discrete Markov process. Second, one can straightforwardly calculate the phase at any particular time through directly performing the minimization in \eqref{minim}. Third, we used this definition of the phase to estimate the time that the system stays in close proximity to the limit cycle, demonstrating that this time is $O\big(\exp(Cb\Omega)\big)$ with probability $O\big(\exp(-Cb\Omega)\big)$, where $b$ is the rate of decay to the limit cycle, and $\Omega$ is the system size. In other words, if $b\Omega \gg 1$, then the system stays close to the limit cycle for extremely long times (much longer than the time over which the diffusion approximation is accurate), with extremely high probability. This analysis also provides a means of assessing the substantial diffusive shifts of the phase that can occur over these very long periods of time. Indeed, the typical order of magnitude of the diffusive shift $\tilde{Y}(t)$ of the phase $\beta$ (defined in \eqref{phase martingale}) is $O\big(\sqrt{\Omega^{-1}t}\big)$. Hence over timescales of $O\big(\exp(Cb\Omega)\big)$, the typical diffusive shift of the phase is $O\big(\Omega^{-1/2}\exp(Cb\Omega/2) \big)$. In most circumstances, one would expect that this shift is much greater than $\omega_0$, the period of the oscillation. In other words, over this timescale one cannot predict the phase $\beta(t)$ with any accuracy at all: all phases are of almost equal likelihood.

 \james{On another note, the probability bound in $\S$5 of this paper differs from the analog for a piecewise-deterministic Markov process in \cite{Bressloff18b}, because the latter does not depend on $b$ to leading order. This is because in the PDMP system, in many circumstances the most likely means by which the system can leave a neighborhood of the limit cycle is to not switch over a significantly long period of time. Since the rate of switching does not depend on the rate of decay towards the limit cycle, the leading order probability bound also does not depend on the rate of decay towards the limit cycle. By contrast, in the CRN network of this paper, if the system does not react over a significant timescale, its distance from the limit cycle does not change.
 
 In general, the methods of this paper are extremely helpful for assessing the robustness of oscillations in chemical reaction networks, which are ubiquitous in cell biology \cite{Bressloff14}. In future work, we plan to extend these methods to chemical reaction networks that support more complicated emergent dynamics.}

\bibliographystyle{siam}

\begin{thebibliography}{10}

\bibitem{Agazzi17} {\sc A. Agazzi, A. Dembo and J. P. Eckmann}, {\em Large deviations theory for Markov jump models of chemical reaction networks.}
Annals of Applied Probability. 28 (2018) Issue 3.

\bibitem{Anderson10} {\sc D. F. Anderson, G. Craciun and T. G. Kurtz,} {\em Product-form stationary distributions for deficiency zero chemical reaction networks.} Bull. Math. Biol. 72 (2010) 1947-1970.


\bibitem{Anderson11} {\sc D. F. Anderson and T. G. Kurtz,} {\em Continuous time Markov chain models for chemical reaction networks.}
In: design and analysis of biomolecular circuits (2011) pp. 3-42.

\bibitem{Anderson15} {\sc D. F. Anderson and T. G. Kurtz}, {\em Stochastic analysis of biochemical systems.} Springer (2015)


\bibitem{Ashwin16}
{\sc P.~Ashwin, S.~Coombes, and R.~Nicks}, {\em Mathematical frameworks for
  oscillatory network dynamics in neuroscience}, J. Math. Neurosci., 6(2)
  (2016).

\bibitem{Boland08} {\sc R.~P. Boland, T.~Galla, and A.~J. McKane,} {\em How limit cycles and quasi-cycles are related in systems with intrinsic noise}. J.. Stat. Mech.: Theory and Experiment, {P09001} (2008) pp. 1-27.


\bibitem{Boland09}
{\sc R.~P. Boland, T.~Galla, and A.~J. McKane}, {\em Limit cycles, complex
  floquet multipliers, and intrinsic noise.}, Phys. Rev. E, 79 (2009),
  p.~051131.
  
  
\bibitem{Bonnin17}
{\sc M.~Bonnin}, {\em Amplitude and phase dynamics of noisy oscillators}, Int.
  J. Circuit Th. Appl., 45 (2017), pp.~636--659.
  
  \bibitem{Bressloff18a} {\sc P. C. Bressloff and J. N. Maclaurin,} {\em  A variational method for analyzing stochastic limit cycle oscillators} {SIAM J. Appl. Dyn. Syst.} 17 (2018) 2205-2233.
  
  \bibitem{Bressloff18b}  {\sc P. C. Bressloff and J. N. Maclaurin,} {\em A variational method for analyzing limit cycle oscillations in stochastic hybrid systems.} {Chaos } 28 (2018) 063105.

\bibitem{Bressloff14} {\sc P. C. Bressloff}, {\em Stochastic processes in cell biology.} Springer (2014)

\bibitem{Bressloff17} {\sc P. C. Bressloff,} {\em Stochastic switching in biology: from genotype to phenotype} J. Phys. A 50 (2017) 133001.

\bibitem{Holmes04}
{\sc E.~Brown, J.~Moehlis, and P.~Holmes}, {\em On the phase reduction and
  response dynamics of neural oscillator populations.}, Neural Comput., 16
  (2004), pp.~673--715.
  
     \bibitem{Elf03}
{\sc J. Elf and M. Ehrenberg,} {\em  Fast evaluation of fluctuations in biochemical networks
  with the linear noise approximation}.
\newblock Genome Res. 13 (2003) pp. 2475-2484.

\bibitem{Epstein98}
{\sc I. R. Epstein and J. A. Pojman,}
{\em An introduction to nonlinear chemical dynamics: oscillations, waves, patterns, and chaos.}
Oxford University Press (1998)

\bibitem{Erm96}
{\sc G.~B. Ermentrout}, {\em Type {I} membranes, phase resetting curves, and
  synchrony}, Neural Comput., 8 (1996), p.~979.



\bibitem{Erm10}
{\sc G.~B. Ermentrout}, {\em Noisy oscillators},
  in Stochastic methods in neuroscience, C~R Laing and G~J Lord, eds., Oxford
  University Press, Oxford (2009).
  


\bibitem{Gardiner09}
{\sc C.~W. Gardiner}, {\em Handbook of stochastic methods, 4th edition} (2009)
  Springer, Berlin.

\bibitem{Glass88}
{\sc L.~Glass and M.~C. Mackey}, {\em From {C}locks to {C}haos} (1988) Princeton Univ
  Press, Princeton.
  
\bibitem{Gillespie77} {\sc D. T. Gillespie,} {\em  Exact stochastic simulation of coupled chemical reactions.} J. Phys. Chem. 81 (1977) 2340-2361.

\bibitem{Gillespie01} {\sc D. T. Gillespie,}  {\em Approximate accelerated stochastic simulation of chemically reacting systems.} J. Chem. Phys. 115 (2001) pp. 1716-1733.



\bibitem{Gillespie13} {\sc D. T. Gillespie, A. Hellander and L. R. Petzold,} {\em Perspective: stochastic algorithms for chemical kinetics.} J. Chem. Phys. {138} (2013) 170901.


\bibitem{Goldobin05}
{\sc D.~S. Goldobin and A.~Pikovsky}, {\em Synchronization and
  desynchronization of self--sustained oscillators by common noise}, Phys. Rev.
  E, 71 (2005), pp.~045201.

\bibitem{Gonze02}
{\sc D.~Gonze, J.~Halloy, and P.~Gaspard}, {\em Biochemical clocks and
  molecular noise: theoretical study of robustness factors.}, J. Chem. Phys.,
  116 (2002), pp.~10997--11010.
  
  \bibitem{Hanggi84}
{\sc P. Hanggi, H. Grabert, P. Talkner and H. Thomas} {\em Bistable systems: master
  equation versus {F}okker--{P}lanck modeling.}
\newblock Phys. Rev. A 29 (1984) pp. 371-378.
  
  
\bibitem{Collins05}
{\sc M. Kaern, T. C. Elston, W. J.  Blake and J. J. Collins,} {\em Stochasticity in gene
  expression: from theories to phenotypes.}
\newblock Nat. Rev. Genetics 6 (2005) pp. 451-464.

\bibitem{keskin18}
{\sc S. Keskin, G.S Devakanmalai, S.B. Kwon, H.T. Vu, Q. Hong, Y.Y. Lee, M. Soltani, A. Singh, A. Ay and E. Ozbudak}
{\em Noise in the vertebrate segmentation clock is boosted by time delays but tamed by notch signaling}
\newblock Cell Reports (2018) 2175 - 2185.

\bibitem{Koeppl11}
{\sc H.~Koeppl, M.~Hafner, A.~Ganguly, and A.~Mehrotra}, {\em Deterministic
  characterization of phase noise in biomolecular oscillators.}, Phys. Biol., 8
  (2011), 055008.
  
  \bibitem{Kumar16}
{\sc R. Kumar and L. Popovic,} {\em  Large deviations for multi-scale jump-diffusion processes}
Stochastic Processes and their Applications (2016)

\bibitem{Kuramoto84}
{\sc Y.~Kuramoto}, {\em Chemical Oscillations, Waves and Turbulence},
  (1984) Springer-Verlag, New-York.
  
\bibitem{Kurtz71} {\sc T. G. Kurtz,} {\em Limit theorems for sequences of jump Markov processes approximating ordinary differential
processes.} Journal of Applied Probability. 8 (1971) pp. 344-356.

\bibitem{Kurtz76}
{\sc T. G. Kurtz,} {\em Limit theorems and diffusion approximations for density dependent
  {M}arkov chains.}
\newblock Math. Prog. Stud. 5 (1976) pp. 67-78.

\bibitem{Kurtz80} {\sc T. G. Kurtz,} {\em Representations of {M}arkov processes as multiparameter changes.} Ann. Prob.8 (1980) pp. 682-715.
  
  \james{
  \bibitem{Maini12} {\sc P.K. Maini and T.E. Woolley and R.E. Baker and E.A. Gaffney and S. Lee}
  {\em Turing's model for biological pattern formation and the robustness problem}
  \newblock Interface focus 2 (4) (2012) pp. 487-496.
  }
   
  \bibitem{Minas17} {\sc G. Minas and D. A. Rand,} {\em Long-time analytic approximation of large stochastic oscillators: simulation, analysis and inference.} PLoS Comput. Biol. 13 (2017) e1005676.

\bibitem{nakao2016}
{\sc H.~Nakao}, {\em Phase reduction approach to synchronization of nonlinear
  oscillators}, Contemporary Physics, 57 (2016), pp.~188--214.

\bibitem{Paulsson05}
{\sc J. Paulsson,} {\em Models of stochastic gene expression.}
\newblock Phys. Life Rev. 2 (2005) pp. 157-175.
\james{
\bibitem{Protter05}
{\sc P. Protter,} {\em Stochastic Integration and Differential Equations. 2nd Edition.}
\newblock Spinger-Verlag Berlin-Heidelberg 2005
}
\bibitem{Raj08}
{\sc A. Raj and A. van Oudenaarden,} {\sc Nature, nurture, or chance: Stochastic gene expression and its consequences.}
\newblock Cell  135 (2008) pp. 216-226.


  
  \bibitem{Suvak12} {\sc O. Suvak and A. Demir,} {\em Phase computations and phase models for discrete molecular oscillators.} J. Bioinf. Syst. Biol.  6 (2012) pp. 1-28.

\bibitem{Swain02}
{\sc M. B. Elowitz, A. J. Levine, E. D. Siggia and P. S. Swain} {\em Stochastic gene
  expression in a single cell.}
\newblock Science 297 (2002)1183-1186.
  
  \bibitem{Tsimiring14} {\sc L. S. Tsimiring,} {\em Noise in biology.} Rep. Prog. Phys. 77 (2014) 026601.
  
 \bibitem{vanKampen92}
{\sc N. G. van Kampen,} {\em Stochastic processes in physics and chemistry.}
\newblock North-Holland, Amsterdam (1992).

\bibitem{Vellela10}
{\sc M. Vellela and H. Qian} Stochastic dynamics and non-equilibrium thermodynamics
  of a bistable chemical system: the {S}chl{\"o}gl model revisited.
\newblock J. R. Soc. Interface 6 (2009) pp. 925-940.


\bibitem{Winfree80}
{\sc A.~Winfree}, {\em The geometry of biological time} (1980) Springer-Verlag, New
  York.


\end{thebibliography}

\end{document}